\documentclass[twoside,10pt]{article}

\usepackage{psfrag}
\usepackage{graphicx}
\usepackage{amssymb}
\usepackage{latexsym}
\usepackage{amsmath}

\newtheorem{theorem}{Theorem}
\newtheorem{proposition}{Proposition}
\newtheorem{definition}{Definition}
\newtheorem{lemma}{Lemma} 
\newtheorem{corollary}{Corollary}
\newtheorem{sublemma}{Sublemma}
\newtheorem{remark}{Remark}

\date{}

\begin{document}

\title{Products of hyperbolic metric spaces}
\maketitle

\begin{center}
{\large Thomas Foertsch* 
\footnote{*supported by SNF Grant 21 - 589 38.99},  \hspace{1cm} Viktor Schroeder
}
 \\
%{\small Institute for Mathematics, University Zurich, E-mail: foertsch@math.unizh.ch}
\end{center}

\begin{abstract}
Let $(X_i,d_i)$, $i=1,2$, be proper geodesic hyperbolic metric spaces.
We give a general construction for a ``hyperbolic product'' $X_1{\times}_hX_2$ which is itself a proper geodesic
hyperbolic metric space and examine its boundary at infinity. 
\end{abstract}

\vspace{0.5cm}

%%%%%%%%%%%%%%%%%%%%%%%%%%%%%%%%%%%%%%%%%%%%%%%%%%%%%%%%%%%%%%%%%%%%%%%%%%%%%%%%%%%%%%%%%%%%%
%%%%%%%%%%%%%%%%%%%%%%%%%%%%%%%%%%%%%%%%%%%%%%%%%%%%%%%%%%%%%%%%%%%%%%%%%%%%%%%%%%%%%%%%%%%%%

\section{Introduction}
\label{sec-intro}

Let $(X_i,d_i)$, $i=1,2$, be proper geodesic hyperbolic metric spaces (for definitions see Section \ref{sec-preliminaries}).
We give a general construction for a ``hyperbolic product'' $X_1{\times}_hX_2$ which is itself a proper geodesic
hyperbolic metric space. \\
This construction only depends on chosen basepoints $z_i\in X_i$ or on points $u_i\in \partial X_i$, where 
$\partial X_i$ is the boundary at infinity. \\

For given points $z_i\in X_i$ consider the set
\begin{displaymath}
Y \; := \; \Big\{ (x_1,x_2)\in X_1\times X_2 \; \Big| \; d_1(x_1,z_1) \; = \; d_2(x_2,z_2)\Big\} \; \subset \; 
X_1 \times X_2
\end{displaymath}
with the induced product metric 
\begin{displaymath}
d_2\Big( (x_1,x_2),(y_1,y_2)\Big) \; = \; \Big( d_1^2(x_1,y_1) \; + \; d_2^2(x_2,y_2){\Big)}^{\frac{1}{2}} .
\end{displaymath}
Let $d:Y\times Y\longrightarrow {\mathbb{R}}^+$ be the corresponding interior metric 
\begin{displaymath}
d(x,y) \; := \; \inf \Big\{ L(c) \; \Big| \; c \; \mbox{connecting} \; x \; \mbox{to} \; y \; \mbox{in} \; Y \Big\} ,
\end{displaymath}
where $L(c)$ denotes the length of $c$ in $(Y,d_e)$. \\
We call $(Y,d)$ the {\it hyperbolic product} of $(X_1,z_1)$ and $(X_2,z_2)$ and denote it also by
\begin{displaymath}
Y \; = \; (X_1,z_1) \; {\times}_h \; (X_2,z_2) . 
\end{displaymath}
This notion is justified by the following
\begin{theorem} \label{theo-main}
Let $X_i$, $i=1,2$, be proper geodesic hyperbolic spaces and $z_i\in X_i$, $i=1,2$. Then 
$Y = (X_1,z_1) {\times}_h (X_2,z_2)$ is also a proper geodesic hyperbolic space and $\partial Y$ is naturally 
homeomorphic to $\partial X_1 \times \partial X_2$.
\end{theorem}

The construction can be carried over in the limit case that the points $z_i$ tend to infinity. This limit case seems to be of
particular interest. \\
Let $u_i\in \partial X_i$ be given. These points give rise to Busemann functions $B_i:X_i\longrightarrow \mathbb{R}$.
Define now
\begin{displaymath}
Y \; := \; \Big\{ (x_1,x_2)\in X_1\times X_2 \; \Big| \; B_1(x_1) \; = \; B_2(x_2) \Big\} \; \subset \; X_1\times X_2
\end{displaymath}
and consider as above the interior metric $d$ on $Y$. \\
We call $(Y,d)$ the hyperbolic product of $(X_1,B_1)$ and $(X_2,B_2)$, denote it by 
\begin{displaymath}
Y \; = \; (X_1,B_1) \; {\times}_h \; (X_2,B_2) ,
\end{displaymath}
and obtain the
\begin{theorem} \label{theo-main2}
Let $X_i$, $i=1,2$, be proper geodesic hyperbolic spaces and $B_i:X_i\longrightarrow \mathbb{R}$
Busemann functions on $X_i$. Then $Y=(X_1,B_1){\times}_h(X_2,B_2)$ is also a proper geodesic hyperbolic space
and $\partial Y$ is naturally homeomorphic to the smashed product $\partial X_1 \wedge \partial X_2$.
\end{theorem}

\begin{remark} \label{rem-1}
\begin{description}
\item[i)] The smashed product $\wedge$ is a standard construction 
for pointed topological spaces. Let $(U_1,u_1)$, $(U_2,u_2)$ be two pointed spaces then the smashed product 
$U_1\wedge U_2$ is defined as $U_1\times U_2/U_1\vee U_2$, where $U_1\times U_2$ 
is the usual product and 
\begin{displaymath}
U_1\vee U_2 \; = \; \Big( \{ u_1\} \times U_2\Big) \; \cup \Big( U_2 \times \{ u_2\} \Big) \; \subset \; U_1\times U_2
\end{displaymath}
is the wedge product canonically embedded in $U_1\times U_2$. Thus $U_1\wedge U_2$ is obtained from $U_1\times U_2$
by collapsing $U_1\vee U_2$ to a point. For example $S^m\wedge S^n=S^{m+n}$.  
\item[ii)] In \cite{brfa} the authors proved that the hyperbolic product of real hyperbolic spaces 
$(\mathbb{H}^{m_1},u_1){\times}_h(\mathbb{H}^{m_2},u_2)$, with $u_i$ in the ideal boundary, 
is isometric to a real hyperbolic space $\mathbb{H}^{m_1+m_2-1}$. \\
In \cite{fs} the authors proved that, more generally, the hyperbolic product $(Y,d)$ of Hadamard manifolds of
pinched negative sectional curvature $-b^2\le K_i \le -a^2 <0$ is hyperbolic. This was done by showing that there 
exists a metric on $Y$ that is bilipschitz to the one induced by the canonical embedding $i:Y\longrightarrow X$ and that 
carries pinched negative sectional curvature. \\
For related results also see \cite{l1}, \cite{l2}, \cite{f1} and \cite{f2}.
\item[iii)] Instead of the Euclidean product metric $d_e$ on $Y$ we could also take e.g. the maximum metric
\begin{displaymath}
d_m\Big( (x_1,x_2), (y_1,y_2)\Big) \; := \; \max \Big\{ d_1(x_1,y_1),d_2(x_2,y_2)\Big\} 
\end{displaymath}
and the corresponding inner metric $d'$ on $Y$. Note that $d_e$ and $d_m$ and therefore also $d$ and $d'$
are bilipschitz related. \\
The following holds in general: If $(Y,d)$ is a proper geodesic metric space and $d'$ is an other interior metric on $Y$ which is
bilipschitz related to $d$, then $(Y,d')$ also is a proper metric space which implies that it is geodesic since 
$(Y,d)$ is a length space (see e.g. Theorem 2.5.23 in \cite{bubui}). In addition $(Y,d)$ is hyperbolic if and only
if $(Y,d')$ is hyperbolic and in that case $\partial (Y,d)$ is homeomorphic to $\partial (Y,d')$. \\
For technical reasons we use in our proof the metric $d_m$ on $X_1\times X_2$.
\item[iv)] Finally note that the hyperbolic product can similar be defined for finitely many factors and 
the analogue of Theorems \ref{theo-main} and \ref{theo-main2} hold in that case.
\end{description}
\end{remark}

{\it Outline of the paper:} \\
In Section \ref{sec-preliminaries} we collect the necessary results on hyperbolic metric spaces. In Sections
\ref{sec-hyp-product} and \ref{sec-boundary} we give a proof of Theorem \ref{theo-main2}. At the end of Section
\ref{sec-hyp-product} we indicate the necessary changes for the situation of Theorem \ref{theo-main}. \\

{\bf Acknowledgment:} We want to thank Urs Lang for useful discussions.

%%%%%%%%%%%%%%%%%%%%%%%%%%%%%%%%%%%%%%%%%%%%%%%%%%%%%%%%%%%%%%%%%%%%%%%%%%%%%%%%%%%%%%%%%%%%%
%%%%%%%%%%%%%%%%%%%%%%%%%%%%%%%%%%%%%%%%%%%%%%%%%%%%%%%%%%%%%%%%%%%%%%%%%%%%%%%%%%%%%%%%%%%%%

\section{Preliminaries}
\label{sec-preliminaries}

%%%%%%%%%%%%%%%%%%%%%%%%%%%%%%%%%%%%%%%%%%%%%%%%%%%%%%%%%%%%%%%%%%%%%%%%%%%%%%%%%%%%%%%%%%%%%

\subsection{Hyperbolicity}
\label{subsec-hyperbolicity}

A metric space $(X,d)$ is called {\it geodesic}, if any two points $x,y\in X$ can be joined by a geodesic segment $\overline{xy}$
that is the image of a geodesic path ${\gamma}_{xy}:[0,d(x,y)]\longrightarrow X$ from $x$ to $y$ which is parameterized by 
arclength. \\
A geodesic metric space is called $\delta${\it -hyperbolic} if for any triangle with geodesic sides in $X$
each side is contained in the $\delta$- neighborhood of the union of the two other sides. \\
The space is called {\it hyperbolic} if it is $\delta$-hyperbolic for some $\delta \ge 0$. \\

Let $X$ be a metric space and $x,y,z\in X$. Then there exist unique $a,b,c\in \mathbb{R}^+_0$ such that 
\begin{displaymath}
d(x,y)=a+b, \hspace{0.5cm} d(x,z)=a+c \hspace{0.5cm} \mbox{and} \hspace{0.5cm} d(y,z)=b+c . 
\end{displaymath}
In fact those numbers are given through
\begin{displaymath}
a \; = \; (y\cdot z)_x \; , \hspace{0.5cm} b \; = \; (x\cdot z)_y \; , \hspace{0.5cm}
\mbox{and} \hspace{0.5cm} c \; = \; (x\cdot y)_z \; , 
\end{displaymath}
where for instance 
\begin{displaymath}
(y\cdot z)_x \; = \; \frac{1}{2} \Big[ d(y,x) \; + \; d(z,x) \; - \; d(y,z) \Big] .
\end{displaymath}
In the case that $X$is geodesic we may consider a geodesic triangle $\overline{xy} \cup \overline{xz} \cup \overline{yz} \subset X$,
where for example $\overline{xy}$ denotes a geodesic segment connecting $x$ to $y$. Given such a triangle we denote
by $\tilde{x}={\gamma}_{xy}(a)$ the unique point on $\overline{yz}$ satisfying $d(\tilde{x},y)=(x\cdot z)_y$ and in the same 
way we define $\tilde{y}\in \overline{xz}$ and $\tilde{z}\in \overline{xy}$. \\
Note that for $X$ being a tree all these points coincide, i.e. $\tilde{x}=\tilde{y}=\tilde{z}$. In general an
upper bound for the distances of these points measures the hyperbolicity of $(X,d)$.

\begin{lemma} \label{lemma-hyp-equiv1}
\begin{description}
\item[i)] If $(X,d)$ is $\delta$-hyperbolic, then 
\begin{displaymath}
d(z,\tilde{z}) \; \le \; c+2\delta , \hspace{1cm} 
d\Big( {\gamma}_{xy}(t),{\gamma}_{xz}(t)\Big) \; \le \; 4\delta \;\;\; \forall t\in [0,a]
\end{displaymath}
and the points $\tilde{x}$, $\tilde{y}$, $\tilde{z}$ have pairwise distance $\le 4\delta$.
\item[ii)] A metric space $(X,d)$ is hyperbolic if and only if there exists a ${\delta}' \in \mathbb{R}_0^+$ such that given any
geodesic triangle $\overline{xy} \cup \overline{xz} \cup \overline{yz} \subset X$ the points $\tilde{x} \in \overline{yz}$, 
$\tilde{y} \in \overline{xz}$ and $\tilde{z} \in \overline{xy}$ as defined above have distance less than ${\delta}'$ to each other.
\end{description}
\end{lemma}
{\bf Proof:} $i)$ By $\delta$-hyperbolicity $d(\tilde{z},\overline{xz})\ge \delta$ or
$d(\tilde{z},\overline{yz})\ge \delta$. By triangle inequality we have in the first case $d(\tilde{z},\tilde{y})\le 2\delta$
and hence $d(\tilde{z},z)\le c+2\delta$. The other case is similar. \\
Assume that there is $t_0\in [0,a]$ with $d({\gamma}_{xy}(t_0),{\gamma}_{xz}(t_0))>4\delta$, then \linebreak
$d({\gamma}_{xy}(t_0-\delta ),{\gamma}_{xz}(t_0-\delta ))>2\delta$ which implies 
$d({\gamma}_{xy}(t_0-\delta ),\overline{xz})>\delta$ and by hyperbolicity 
$d({\gamma}_{xy}(t_0-\delta ),\overline{yz})<\delta$. 
Let $p\in \overline{yz}$ be a point of minimal distance to $\gamma (t_0-\delta)$. By triangle inequality 
$d(p,y)\ge b$ and hence $d(p,z)\le c$. Thus 
\begin{eqnarray*}
a \; + \; c \; = \; d(x,z) & \le & 
(t_0-\delta ) \; + \; d\Big( {\gamma}_{xy}(t_0-\delta ), p\Big) \; + \; d(p,z) \\
& < & a \; + \; c \; ;
\end{eqnarray*}
a contradiction. Since the corresponding estimate holds for the other sides as well the points $\tilde{x}$,
$\tilde{y}$ and $\tilde{z}$ have pairwise distance $\le 4\delta$. \\
For $i)$ compare to Proposition III.1.17 in \cite{brihae}.
\hfill{$\Box$}

%%%%%%%%%%%%%%%%%%%%%%%%%%%%%%%%%%%%%%%%%%%%%%%%%%%%%%%%%%%%%%%%%%%%%%%%%%%%%%%%%%%%%%%%%%%%%

\subsection{$T$-functions}
\label{subsec-t-functions}

\begin{definition}
Let $\alpha ,\omega \in \mathbb{R}$ and $I:=[\alpha ,\omega ]$.
\begin{description} 
\item[i)] A function $f:I\longrightarrow \mathbb{R}$ is called a $T$-function if $f$ is continuous and there exists 
$\alpha +a\in [\alpha ,\omega ]$ such that the restrictions $f|_{(\alpha ,\alpha +a)}$ and $f|_{(\alpha +a,\omega )}$ 
are differentiable with 
\begin{displaymath}
f'|_{(\alpha ,\alpha +a)} \; \equiv \; -1 \hspace{1cm} \mbox{and} \hspace{1cm} f'|_{(\alpha +a,\omega )} \; \equiv \; 1 .
\end{displaymath}
\item[ii)] A function $f:I\longrightarrow \mathbb{R}$ is called a $\delta$-$T$-function, $\delta \in \mathbb{R}_0^+$, 
if there exists a $T$-function $g:I\longrightarrow \mathbb{R}$ such that $||f-g||_{\sup} \, < \, \delta$.
\item[iii)] Let $X$ be a geodesic metric space. A function $f:X\longrightarrow \mathbb{R}$ is called a
$\delta$-$T$-function, if for any geodesic segment $\gamma :[\alpha ,\omega ]\longrightarrow \mathbb{R}$ 
the function $f\circ \gamma$ is a 
$\delta$-$T$-function.
\end{description} 
\end{definition}

\begin{remark} It is straight forward to check that
\begin{description} 
\item[i)] Every $T$-function $f$ is convex and Lipschitz with Lipschitz constant $l=1$.
\item[ii)] For $[\alpha ,\omega ]\subset \mathbb{R}$ and $t_1,t_2\in \mathbb{R}$ with $|t_1-t_2|\le |\alpha -\omega |$ 
there exists a unique $T$-function with $f(\alpha )=t_1$ and $f(\omega )=t_2$. Indeed there are unique 
$a,b,c\in \mathbb{R}_0^+$ such that $a+b=|\alpha -\omega |$, $a+c=t_1$, $b+c=t_2$. These are given via
\begin{eqnarray*}
a & = & \frac{1}{2} \Big( |\alpha -\omega | \; + \; t_1 \; - \; t_2\Big) , \\
b & = & \frac{1}{2} \Big( |\alpha -\omega | \; + \; t_2 \; - \; t_1\Big) , \\
c & = & \frac{1}{2} \Big( t_1 \; + \; t_2 \; - \; |\alpha -\omega |\Big) . 
\end{eqnarray*}
Now $f$ satisfies
\begin{displaymath}
f(\alpha ) \; = \; t_1 , \hspace{0.5cm} f(\alpha +a) \; = \; c , \hspace{0.5cm} f(\omega ) \; = \; t_2 .
\end{displaymath}
\item[iii)] $f$ is a ($\delta$-)$T$-function $\Longrightarrow$ $f+const$ is a  ($\delta$-)$T$-function.
\item[iv)] A limit of a sequence of ($\delta$-)$T$-functions is a ($\delta$-)$T$-function.
\end{description}
\end{remark}

\begin{lemma}  \label{lemma-hyp-equiv2}
Let $X$ be a geodesic metric space. Then the following are equivalent:
\begin{description}
\item[i)] $X$ is hyperbolic.
\item[ii)] There exists a $\delta \in \mathbb{R}^+$ such that
for all $x\in X$ the function 
\begin{displaymath}
d_x:  X \longrightarrow  \mathbb{R} \, , \hspace{1cm} y\longmapsto d(x,y)
\end{displaymath}
is a $\delta$-$T$-function.
\end{description}
\end{lemma}
{\bf Proof:} ``$\Longrightarrow$'' Let $X$ be $\delta$-hyperbolic, $x,y,z\in X$, $d(x,z)=a+c$, $d(y,z)=b+c$ and 
$\gamma :[0,a+b]\longrightarrow \overline{xy}$ an arc length
parameterized geodesic connecting $x$ to $y$. \\
Now consider the $T$-function $f:[0,a+b]\longrightarrow \mathbb{R}^+$ determined by 
\begin{displaymath} 
f(0)=d(z,\gamma (0))=a+c \hspace{0.5cm} \mbox{and} \hspace{0.5cm}
f(a+b)=d(z,\gamma (a+b))=b+c. 
\end{displaymath}
Note that $f(a)=c$ and $(d_z\circ \gamma )(a)\le c+4\delta$ by Lemma \ref{lemma-hyp-equiv1} $ii)$.

From the fact that 
$d_z\circ \gamma : [0,a+b]\longrightarrow \mathbb{R}^+$ is 1-Lipschitz it immediately follows that
\begin{displaymath}
||d_z\circ \gamma \; - \; f ||_{\sup} \; \le \; 4 \delta .
\end{displaymath}
Hence $d_z\circ \gamma$ is a $4\delta$-$T$-function. \\
``$\Longleftarrow$'' Let now $X$ satisfy condition $ii)$. We show that $X$ is hyperbolic using 
the criterion of Lemma \ref{lemma-hyp-equiv1} $i)$: \\
For $x,y,z\in X$ and geodesic segments $\overline{xy}$, $\overline{xz}$ and $\overline{yz}$ connecting these points, condition
$ii)$ yields $c\le d(z,\tilde{z})\le c+\delta$. We now consider the geodesic triangle 
$\overline{x\tilde{z}} \cup \overline{z\tilde{z}} \cup \overline{xz}$, where $\overline{x\tilde{z}}\subset \overline{xy}$ and
$\overline{z\tilde{z}}$ is any geodesic segment connecting $z$ to $\tilde{z}$. For $\hat{z}\in \overline{xz}$
satisfying $d(z,\hat{z})=(x\cdot \tilde{z})_z$ condition $ii)$ gives $d(\tilde{z}, \hat{z})< \delta$.
Furthermore one has $c<d(z,\hat{z})<c+\frac{\delta}{2}$ and therefore $d(\tilde{y},\hat{z})<\frac{\delta}{2}$.
Thus we achieve
\begin{displaymath}
d(\tilde{z},\tilde{y}) \; \le \; d(\tilde{z},\hat{z}) \; + \; d(\hat{z},\tilde{y}) \; \le \;
\delta \; + \; \frac{\delta}{2} \; = \; \frac{3}{2} \delta .
\end{displaymath}
The same argument of course yields $d(\tilde{z},\tilde{x})\le \frac{3}{2} \delta$ and 
$d(\tilde{y},\tilde{x})\le \frac{3}{2}\delta$. By Lemma \ref{lemma-hyp-equiv1} $ii)$ we obtain the result.  \hfill $\Box$

\begin{figure}[htbp]
\psfrag{0}{$0$}
\psfrag{a}{$a$}
\psfrag{a+b}{$a+b$}
\psfrag{c}{$c$}
\psfrag{a+c}{$a+c$}
\psfrag{b+c}{$b+c$}
\psfrag{delta}{$\delta$}
\psfrag{f}{$f$}
\psfrag{f+delta}{$f+\delta$}
\psfrag{d_zcircgamma}{$d_z \circ \gamma$}
\includegraphics[width=0.6\columnwidth]{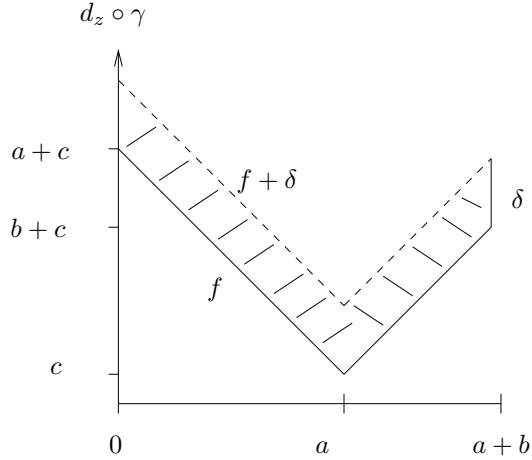}
\caption{The graph of $d_z\circ \gamma$
lies in the region between the graphs of $f$ and $f+\delta$.} 
\end{figure}

%%%%%%%%%%%%%%%%%%%%%%%%%%%%%%%%%%%%%%%%%%%%%%%%%%%%%%%%%%%%%%%%%%%%%%%%%%%%%%%%%%%%%%%%%%%%%

\subsection{The boundary at infinity and Busemann functions}
\label{subsec-b-functions}

We associate to a hyperbolic space $(X,d)$ a boundary $\partial X$ at infinity. There are different descriptions in the
literature (see e.g. \cite{beka}) all of which coincide for proper geodesic hyperbolic spaces. \\
We choose a basepoint $z\in X$. We say that a sequence $\{ x^k{\}}_{k\in \mathbb{N}}$ of points in $X$ 
{\it converges to infinity}, if 
\begin{displaymath}
\liminf\limits_{k,l\longrightarrow \infty} \; (x^k\cdot x^l)_z \; = \; \infty .
\end{displaymath}
Two sequences  $\{ x^k{\}}_{k\in \mathbb{N}}$ and  $\{ y^k{\}}_{k\in \mathbb{N}}$ converging to infinity
are equivalent, \linebreak  $\{ x^k{\}}_{k\in \mathbb{N}} \sim \{ y^k{\}}_{k\in \mathbb{N}}$ if
\begin{displaymath}
\liminf\limits_{k,l\longrightarrow \infty} \; (x^k\cdot y^l)_z \; = \; \infty .
\end{displaymath}
One shows that $\sim$ is an equivalence relation and defines $\partial X$ as the set of equivalence classes.
We write $[\{ x^k\} ]\in \partial X$ for the corresponding class. One can also show that for every $x\in X$
and $v\in \partial X$ there is a geodesic ray ${\gamma}_{xv}:[0,\infty] \longrightarrow X$ parameterized by arclength
with ${\gamma}_{xv}(0)=x$ and $[\{ {\gamma}_{xv}(k)\} ]=v$. \\

For $v\in \partial X$ and $r>0$ one defines
\begin{displaymath}
U(v,r) \; := \; 
\Big\{ {\textstyle
w\in \partial X \; \Big| \; \exists \{ x^k\} , \{ y^k\} \; \mbox{s.t.} \;
[\{ x^k\} ]=v, \;  [\{ y^k\} ]=w, \; \liminf\limits_{k,l\longrightarrow \infty} \; (x^k\cdot y^l)_z \; \ge \; r }
\Big\} . 
\end{displaymath} 
On $\partial X$ we consider the topology generated by $U(v,r)$, $v\in \partial X$, $r>0$. \\

If $X$ is proper, $\partial X$ is a compact topological space. We now also fix a basepoint $u\in \partial X$ and
a geodesic ray ${\gamma}_{zu}$ from $z$ to $u$. The function $B:X\longrightarrow \mathbb{R}$,
$B(x):=\lim_{t\rightarrow \infty}[d(x,{\gamma}_{zu}(t))-t]$ is called the {\it Busemann function associated to} ${\gamma}_{zu}$.
$B$ is a $\delta$-$T$-function as a limit of $\delta$-$T$-functions.

\begin{definition}
A geodesic ray $\gamma :[0,\infty ) \longrightarrow X$ is called a $B$-ray if and only if $\gamma$ is parameterized by arc 
length and 
\begin{displaymath}
B\Big( \gamma (t)\Big) \; = \; B\Big( \gamma (0)\Big) \; - \; t \hspace{0.5cm} \forall t\in [0,\infty ).
\end{displaymath}
\end{definition}

By a standard limit argument we obtain the
\begin{lemma} \label{lemma-b-ray-existence}
Let $X$ be a complete, locally compact, hyperbolic metric space and $B$ a Busemann function on $X$. Then for every $x\in X$
there exists a $B$-ray $\gamma$ with $\gamma (0)=x$.
\end{lemma} 

Let now $v\in \partial X\setminus \{ u\}$ and consider $B$-rays ${\gamma}_{{\gamma}_{zv}(t)}$ starting at
${\gamma}_{zv}(t)$. Then these rays subconverge as sets to a geodesic from $v$ to $u$ and by suitable reparameterization 
we obtain the existence of a geodesic ${\gamma}_{vu}:\mathbb{R} \longrightarrow X$ with
\begin{displaymath}
B\Big( {\gamma}_{vu}(t)\Big) \; = \; -t, \hspace{0.5cm}
\Big[ \{ {\gamma}_{vu}(k) \} \Big] \; = \; u \hspace{0.5cm} \mbox{and} \hspace{0.5cm}
\Big[ \{ {\gamma}_{vu}(-k) \} \Big] \; = \; v .
\end{displaymath}
The hyperbolicity of $X$ also implies that there exists a constant $C$ (depending on $u$, $v$, $z$) such that
\begin{displaymath}
d\Big( {\gamma}_{z}(t),{\gamma}_{vu}(t)\Big) \; \le \; C \hspace{0.5cm} \mbox{and} \hspace{0.5cm}
d\Big( {\gamma}_{zv}(t),{\gamma}_{vu}(-t)\Big) \; \le \; C \hspace{0.5cm} \forall t\ge 0 .
\end{displaymath}
Let now $x,y\in X$ and ${\gamma}_{x}, {\gamma}_{y}$ be $B$-rays starting at $x,y$. Then
${\gamma}_{x}, {\gamma}_{y}$ and a geodesic segment $\overline{xy}$ form an ideal triangle with vertices
$x,y$ and $u={\gamma}_{zu}(\infty)\in \partial X$. We want to look for points $\tilde{u}, \tilde{x},\tilde{y}$ as for finite 
triangles. Clearly there are $a,b\ge 0$ such that $d(x,y)=a+b$ and $B({\gamma}_{x}(a))=B({\gamma}_{x}(b))$. Indeed
\begin{eqnarray*}
a \; = \; (y\cdot B)_x & := & \frac{1}{2} \Big( d(y,x) \; + \; B(x) \; - \; B(y) \Big) \;\;\; \mbox{and} \\
b \; = \; (x\cdot B)_y & := & \frac{1}{2} \Big( d(x,y) \; + \; B(y) \; - \; B(x) \Big) .
\end{eqnarray*} 
Let $\tilde{u}\in \overline{xy}$ be the point with $d(x,\tilde{u})=a$ and let $\tilde{y}={\gamma}_x(a)$,
$\tilde{x}={\gamma}_y(b)$.

\begin{figure}[htbp]
\psfrag{a}{$a$}
\psfrag{b}{$b$}
\psfrag{c}{$c$}
\psfrag{x}{$x$}
\psfrag{y}{$y$}
\psfrag{z}{$z$}
\psfrag{u}{$u$}
\psfrag{tildex}{$\tilde{x}$}
\psfrag{tildey}{$\tilde{y}$}
\psfrag{tildez}{$\tilde{z}$}
\psfrag{tildeu}{$\tilde{u}$}
\includegraphics[width=0.9\columnwidth]{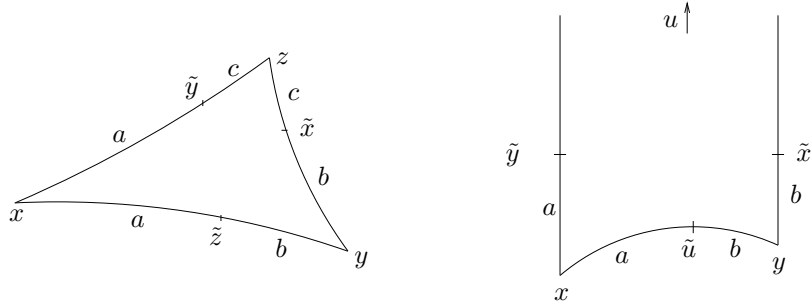}
\caption{Triangle and ideal triangle in the hyperbolic space.}
\end{figure}

\begin{lemma} \label{lemma-xyI}
\begin{description}
\item[i)] The points $\tilde{x}$, $\tilde{y}$ and $\tilde{u}$ have pairwise distance $\le 8\delta$.
\item[ii)] For all $t\ge 0$ it holds $d({\gamma}_{x}(a+t),{\gamma}_y(b+t))\le 8\delta$.
\end{description}
\end{lemma}
{\bf Proof:} Let $B(x)=\lim_{t\rightarrow \infty}(d(x,{\gamma}_{zu}(t))-t)$, $z_i={\gamma}_{zu}(i)$, $i\in \mathbb{N}$,
and consider the triangle $x,y,z_i$ with corresponding values $a_i,b_i,c_i\in \mathbb{R}^+$ as well as points
\begin{displaymath}
\hat{z}_i \; = \; {\gamma}_{xy}(a_i), \hspace{0.5cm} \hat{y}_i \; = \;  {\gamma}_{xz_i}(a_i) 
\hspace{0.5cm} \mbox{and} \hspace{0.5cm}  \hat{x}_i \; = \;  {\gamma}_{yz_i}(b_i) . 
\end{displaymath}
Clearly one has $a_i\longrightarrow a$ and $b_i\longrightarrow b$. \\
Consider also the triangles $x,{\gamma}_x(a_i),z_i$. Note that 
$|d(z_i,{\gamma}_x(a_i))-d(z_i,\hat{y}_i)|\longrightarrow 0$ which implies by Lemma \ref{lemma-hyp-equiv1} $ii)$
\begin{displaymath}
\limsup\limits_{i\longrightarrow \infty} d\Big( {\gamma}_x(a_i), \tilde{y}_i\Big) \; \le \; 2\delta .
\end{displaymath}
In the same way we obtain
\begin{displaymath}
\limsup\limits_{i\longrightarrow \infty} d\Big( {\gamma}_y(b_i), \tilde{x}_i\Big) \; \le \; 2\delta .
\end{displaymath}
Since $d(\tilde{x}_i,\tilde{y}_i)\le 4\delta$ by Lemma \ref{lemma-hyp-equiv1} we obtain $i)$. \\
The proof of $ii)$ is similar.
\hfill $\Box$ \\

We need the following

\begin{lemma} \label{lemma-sigma}
Let $x,y\in X$ and $\sigma :[0,d(x,y)]\longrightarrow X$ a curve parameterized by arclength such that
$\sigma (0)=x$ and $d(\sigma (d(x,y)),y)\le R$. Then
\begin{displaymath}
d\Big( {\gamma}_{xy}(t), \sigma (t)\Big) \; \le \; \frac{3}{2} R \; + \; 4\delta \;\;\; \forall t\in [0,d(x,y)].
\end{displaymath}
\end{lemma}
{\bf Proof:} For fixed $t\in [0,d(x,y)]$ consider the triangle $x,y,z:=\sigma (t)$. Since $d(x,z)\le t$ and
$d(z,y)\le d(x,y)-t+R$ we have $c=(x\cdot y)_z\le \frac{R}{2}$ and hence there exists $\tilde{z}={\gamma}_{xy}(t')$
with $d(z,\tilde{z})\le \frac{R}{2}+4\delta$. Note that $t'=d(x,z)-c$ and $d(x,y)-t'=d(z,y)-c$. Thus
$|t-t'|\le R$ from which the claim follows.
\hfill $\Box$

%%%%%%%%%%%%%%%%%%%%%%%%%%%%%%%%%%%%%%%%%%%%%%%%%%%%%%%%%%%%%%%%%%%%%%%%%%%%%%%%%%%%%%%%%%%%%%

\subsection{A Morse estimate}
\label{subsec-morse}

We need an estimate whose proof is similar to the proof of the Morse inequality.

\begin{lemma} \label{lemma-morse-estimate}
Let $(X,d)$ be $\delta$-hyperbolic, $x,y\in X$ and $\gamma :[0,1]\longrightarrow X$ be a continuous path from $x$ to $y$.
If there exists a point $p={\gamma}_{xy}(s_0)\in \overline{xy}$ such that $d(p,\gamma (t))\ge R$ for all
$t\in [0,1]$ and $R>90\delta $, then 
\begin{displaymath}
L(\gamma ) \; \ge \; d(x,y) \; + \; \frac{1}{20\delta }R^2 .
\end{displaymath}
\end{lemma}
{\bf Proof:} Define $a(t):=(y\cdot \gamma (t))_x\in [0,d(x,y)]$. Since $a$ is continuous, 
$a(0)=0$ and $a(1)=d(x,y)$ there are $0<t_-<t_+<1$ such that 
\begin{displaymath}
a(t_-) \; = \; s_0 \; - \; \frac{R}{2} \; , \hspace{0.5cm} a(t_+) \; = \; s_0 \; + \; \frac{R}{2} \; , 
\hspace{0.5cm}  a(t)\in \Big[ s_0-\frac{R}{2},s_0+\frac{R}{2}\Big] \;\;\; \forall t_-<t<t_+ .
\end{displaymath}
Choose $k+1:=[\frac{R}{8\delta}]+1$ points $s_1\le ...\le s_{k+1} \in [s_0-\frac{R}{2},s_0+\frac{R}{2}]$ 
such that $|s_{i+1}-s_i|\ge 12\delta$ and let $t_i\in [t_-,t_+]$ be points $t_1\le t_2\le ...\le t_{k+1}$ with
$a(t_i)=s_i$. Now
\begin{displaymath}
L(\gamma ) \; \ge \; d\Big( x,{\gamma}(t_-)\Big) \; + \; \sum\limits_{i=1}^k \, d\Big( {\gamma}(t_{i+1}),{\gamma}(t_i)\Big)
\; + \; d\Big( {\gamma}(t_+),y\Big) .
\end{displaymath}
By construction
\begin{displaymath}
d\Big( x,\gamma (t_-)\Big) \ge a(t_-) \hspace{0.5cm} \mbox{and} \hspace{0.5cm}
d\Big( y,\gamma (t_+)\Big) \ge d(x,y) \; - \; a(t_-) . 
\end{displaymath}
Thus $d(x,{\gamma}(t_-))+d({\gamma}(t_+),y)\ge d(x,y)-R$. \\
Let $q_i\in \overline{xy}$ be a point such that $d(\gamma (t_i),q_i)=d(\gamma (t_i),\overline{xy})$.
Since $d_{\gamma (t_i)}\circ {\gamma}_{xy}$ is a $4\delta$-$T$-function by Lemma \ref{lemma-hyp-equiv2}, a minimum of
this function is assumed in distance $\le 4\delta$ of the corresponding $T$-function. 
Thus $d(q_i,{\gamma}_{xy}(s_i))\le 4\delta$, which implies $d(q_i,q_{i+1})\ge 4\delta$ since 
$d({\gamma}_{x,y}(s_{i+1}),{\gamma}_{xy}(s_i))\ge 12\delta$. \\
By Lemma 8.4.23 in \cite{bubui} we obtain
\begin{displaymath}
d\Big( {\gamma}(t_{i+1},{\gamma}(t_i)) \; \ge \; 
d\Big( {\gamma}(t_{i+1}),q_{i+1}\Big) \; + \; d\Big( {\gamma}(t_i),q_i\Big) \; - \; 4\delta .
\end{displaymath}
Since $d(\gamma (t_i),p)\ge R$ and $d(q_i,p)\le \frac{R}{2}+4\delta$ we obtain
$d(\gamma (t_{i+1}),\gamma (t_i)) \; \ge \; R-12\delta$. Thus
\begin{eqnarray*}
L(\gamma ) & \ge & d(x,y) \; + \; \Big[ \frac{R}{12\delta}\Big] (R \; - \; 12\delta ) \; - \; R \\
& \ge & d(x,y) \; + \; \frac{1}{20\delta} R^2  
\end{eqnarray*}
for $R\ge 90 \delta$.
\hfill $\Box$

\begin{figure}[htbp]
\psfrag{qi}{${\scriptstyle q_i}$}
\psfrag{qi+1}{${\scriptstyle q_{i+1}}$}
\psfrag{ge4d}{${\scriptstyle \ge 4\delta}$}
\psfrag{geR/2+4d}{${\scriptstyle \ge \frac{R}{2}+4\delta}$}
\psfrag{R/2+4dle}{${\scriptstyle \frac{R}{2}+4\delta \le}$}
\psfrag{geR-12d}{${\scriptstyle \ge R-12\delta}$}
\psfrag{g(ti)}{${\scriptstyle \gamma (t_i)}$}
\psfrag{g(ti+1)}{${\scriptstyle \gamma (t_{i-1})}$}
\psfrag{gxy(s0+R/2)}{${\scriptstyle {\gamma}_{xy}(s_0+\frac{R}{2})}$}
\psfrag{gxy(s0-R/2)}{${\scriptstyle {\gamma}_{xy}(s_0-\frac{R}{2})}$}
\psfrag{p=gxy(s0)}{${\scriptstyle p={\gamma}_{xy}(s_0)}$}
\includegraphics[width=0.9\columnwidth]{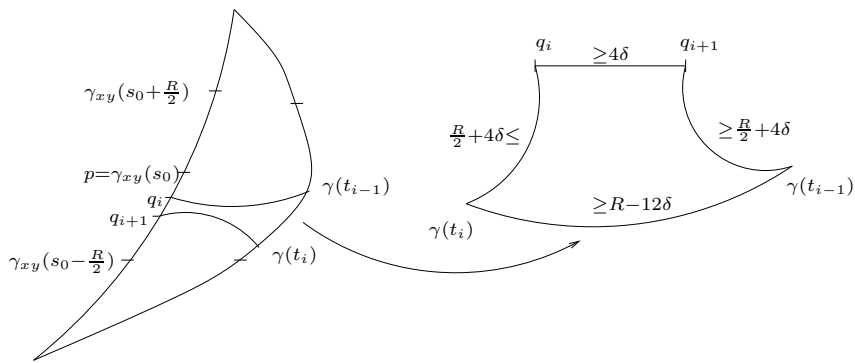}
\caption{This figure shows how Lemma  8.4.23 in \cite{bubui} is applied to the segments
${\gamma}_{xy}|_{[s_i,s_{i+1}]}$ and ${\gamma}|_{t_i,t_{i+1}}$ in the proof of Lemma \ref{lemma-morse-estimate}.} 
\end{figure}

%%%%%%%%%%%%%%%%%%%%%%%%%%%%%%%%%%%%%%%%%%%%%%%%%%%%%%%%%%%%%%%%%%%%%%%%%%%%%%%%%%%%%%%%%%%%%
%%%%%%%%%%%%%%%%%%%%%%%%%%%%%%%%%%%%%%%%%%%%%%%%%%%%%%%%%%%%%%%%%%%%%%%%%%%%%%%%%%%%%%%%%%%%%

\section{The hyperbolic product}
\label{sec-hyp-product}

In this section we prove the first part of Theorem \ref{theo-main2}, which is equal to Proposition \ref{prop-yhyperbolic}. \\

Let $(X_i,d_i)$ be ${\delta}_i$-hyperbolic spaces, $i=1,2$, and $\delta := \max \{ {\delta}_1,{\delta}_2\}$.
Let further $B_i:X_i\longrightarrow \mathbb{R}$ be Busemann functions on $X_i$. We study the set
\begin{displaymath}
Y \; := \; \Big\{ (x_1,x_2)\in X_1\times X_2 \; \Big| \; B_1(x_1) \; = \; B_2(x_2)\Big\} .
\end{displaymath}
On $Y$ we consider the maximum metric $d_m:Y\times Y \longrightarrow \mathbb{R}$,
\begin{displaymath}
d_m \Big( (x_1,x_2) , (x_1',x_2')\Big) \; := \; \max \Big\{ d_1(x_1,x_1'), d_2(x_2,x_2')\Big\} .
\end{displaymath}
For our use of $d_m$ instead of $d_e$ compare Remark \ref{rem-1} $iii)$ in the Introduction. \\
Let $p,p'\in Y$. We first construct two curves in $Y$ between $p$ and $p'$, the {\it $\Gamma$-curve} ${\Gamma}_{pp'}$
and the {\it continuous $\Gamma$-curve} ${\Gamma}^c_{pp'}$. \\
The advantage of ${\Gamma}_{pp'}$ is that this curve is conceptually easy to understand. However, ${\Gamma}_{pp'}$ is
in general not continuous, ${\Gamma}^c_{pp'}$ is a continuous variation of ${\Gamma}_{pp'}$. \\
Let $p=(p_1,p_2)$, $p'=(p_1',p_2')$ and ${\gamma}_i,{\gamma}_i'$ be $B_i$-rays starting at $p_i,p_i'$. Let further
$\gamma ,{\gamma}' :\mathbb{R}_0^+\longrightarrow Y$ be the geodesic rays
\begin{displaymath}
\gamma (t) \; = \; \Big( {\gamma}_1(t),{\gamma}_2(t)\Big) \; , \hspace{1cm} {\gamma}'(t) \; = \; 
\Big( {\gamma}_1'(t),{\gamma_2'(t)}\Big) .
\end{displaymath}
These geodesics are parameterized by constant speed 1. We set
\begin{eqnarray*}
a_i \; := \; (p_i'\cdot B_i)_{p_i} & = & 
\frac{1}{2} \Big( d(p_i,p_i') \; + \; B_i(p_i) \; - \; B_i(p_i')\Big) \;\;\;\; \mbox{and} \\
b_i \; := \; (p_i\cdot B_i)_{p_i'} & = & 
\frac{1}{2} \Big( d(p_i,p_i') \; + \; B_i(p_i') \; - \; B_i(p_i)\Big) \; , 
\end{eqnarray*}
such that $a_i+b_i=d_i(p_i,p_i')$. Let $a:=\max \{ a_1,a_2\}$ and $b:= \max \{b_1,b_2\}$, then $a+b=d_m(p,p')$.
We define by a slight abuse of notation
\begin{displaymath}
\Gamma \; = \; {\Gamma}_{pp'} \; := \; \gamma |_{[0,a]} \; * \; {{\gamma}' |_{[0,b]}}^{-1} .
\end{displaymath}
Note that ${\Gamma}_{pp'}$ is not necessarily continuous, since ${\gamma}(a)=({\gamma}_1(a),{\gamma}_2(a))$ is not
necessarily equal to  ${\gamma}'(b)=({\gamma}'_1(b),{\gamma}'_2(b))$. However, $d_i({\gamma}_i(a),{\gamma}_i'(b))\le 8{\delta}_i$
by Lemma \ref{lemma-xyI}. \\

The curve ${\Gamma}^c_{pp'}$ is a continuous modification of ${\Gamma}_{pp'}$ and defined as 
\begin{displaymath}
{\Gamma}^c_{pp'} \; := \; {\gamma}|_{[0,a+2\delta ]} \; * \; {\Gamma}_1 \; * \; {\Gamma}_2 \; * \;
{{\gamma}'|_{[0,b+2\delta]}}^{-1} ,
\end{displaymath}
where ${\Gamma}_1$ is a continuous curve in $Y$ form $({\gamma}_1(a+2\delta ),{\gamma}_2(a+2\delta ))$ to
$({\gamma}_1'(b+2\delta ),{\gamma}_2(a+2\delta ))$ and ${\Gamma}_2$ a continuous curve in $Y$ from
$({\gamma}_1'(b+2\delta ),{\gamma}_2(a+2\delta ))$ to $({\gamma}_1'(b+2\delta ),{\gamma}_2'(b+2\delta ))$ 
given in the following way: \\
Let ${\eta}_1:[{\alpha}_1,{\beta}_1]\longrightarrow X_1$ be a geodesic from ${\gamma}_1(a+2\delta )$ to
${\gamma}_1'(b+2\delta )$. Note that $L({\eta}_1)\le 8\delta$ and $B_1({\eta}_1({\alpha}_1))=B_2({\eta}_1({\alpha}_2))$.
Since $B_1$ is 1-Lipschitz and a $4\delta$-$T$-function we obtain
\begin{displaymath}
B_1\Big( {\eta}_1(t)\Big) \; \le \; B_1\Big( {\eta}_1({\alpha}_1)\Big) \; + \; 2\delta \; = \; B_1(p_1) \; - \; a .
\end{displaymath}
Thus $-B_1({\eta}_1(t))+B_1(p_1)\ge a\ge 0$ and
\begin{displaymath}
{\Gamma}_1(t) \; = \; \Big( {\eta}_1(t), {\gamma}_2 \Big( -B_1({\eta}_1(t))+B_1(p_1)\Big) \Big)
\end{displaymath}
is well defined. By construction ${\Gamma}_1(t)\in Y$ and $L_{d_m}({\Gamma}_1)\le 8\delta$. \\
In a similar way one constructs ${\Gamma}_2$. \\

We can easily estimate the length of ${\Gamma}^c_{pp'}$ and obtain the
\begin{lemma} \label{lemma-estimate-Gammac}
Given two points $p,p'\in Y$ the continuous curve ${\Gamma}^c={\Gamma}^c_{pp'}$ has length
\begin{displaymath}
L({\Gamma}^c) \; \le \; d_m(p,p') \; + \; 20 \delta .
\end{displaymath}
\end{lemma}

This immediately implies the

\begin{proposition} \label{prop-d-dm}
Given two points $p,p'\in Y$, it holds
\begin{displaymath}
d_m(p,p') \; \le \; d(p,p') \; \le \; d_m(p,p') \; + \; 20 \delta .
\end{displaymath}
\end{proposition}

Thus $(Y,d_m)$ and $(Y,d)$ are quasi-isometric and hence bilipschitz on a large scale. But also on a local scale they induce
the same topology:

\begin{lemma} \label{lemma-same-topology}
The metrics $d$ and $d_m|_Y$ induce the same topology on $Y$.
\end{lemma}
{\bf Proof:} We need to show that
\begin{displaymath}
\lim\limits_{i\longrightarrow \infty} \, d_m|_Y (y_i,y_0) \; = \; 0 \hspace{0.5cm}
\Longrightarrow \hspace{0.5cm}
\lim\limits_{i\longrightarrow \infty} \, d(y_i,y_0) \; = \; 0 .
\end{displaymath}
Thus it suffices to prove the
\begin{sublemma} \label{sublemma1}
Let $\{ y_i{\}}_{i\in \mathbb{N}}$ be a sequence in $Y$ that converges to $y_0\in Y$ with respect to 
$d_e|_Y$. Then for all $\epsilon >0$ there exists $\rho (\epsilon )>0$ such that for all $y_i\in Y$ satisfying 
$d_e|_Y(y_i,y_0)<\rho$ there exists a curve ${\Gamma}_i$ in $Y$ connecting $y_i$ to $y_0$ of length $L({\Gamma}_i )\le \epsilon$.
\end{sublemma}
A sequence $\{ y_i{\}}_{i\in \mathbb{N}} = \{ (y_{i1},y_{i2}){\}}_{i\in \mathbb{N}}$ in $Y\subset X_1\times X_2$ converges
with respect to $d_m|_Y$ if and only if the sequences $\{ y_{ij}{\}}_{i\in \mathbb{N}}$ in $X_j$ converge with respect to
$d_j$, $j=1,2$. \\
Define 
\begin{displaymath}
K_j \; := \; 
\Big\{ 
\gamma (\frac{\epsilon}{4}) \; \Big| \; \gamma \; \mbox{is a} \; B_j\mbox{-ray with} \; \gamma (0) \; = \; y_{0j}
\Big\}
\end{displaymath}
and let ${\gamma}_{ij}$ be $B_j$-rays with ${\gamma}_{ij}(0)=y_{ij}$. The local compactness of $X_j$ implies the \\
{\it Claim:} There exists $\rho >0$ such that $d_j(y_{ij},y_{0j})<\delta$ implies 
$d({\gamma}_{ij}(\frac{\epsilon}{4}),K_j)\le \frac{\epsilon}{8}$. \\

Let now $d_m(y_i,y_0)<\rho$, then $d_j(y_{ij},y_{0j})<\rho$ and thus by the claim there are $B_j$-rays ${\gamma}_j$
starting at $y_{0j}$ such that $d({\gamma}_{ij}(\frac{\epsilon}{4}),{\gamma}_j(\frac{\epsilon}{4}))<\frac{\epsilon}{8}$.
Similar to the construction of ${\Gamma}^c_{pp'}$ one now finds a continuous path 
${\tilde{\Gamma}}^c$ in $Y$ connecting successively the points
\begin{displaymath}
\Big( y_{i1},y_{i2}\Big) , \;\; 
\Big( {\gamma}_{i1}(\frac{\epsilon}{8}), {\gamma}_{i2}(\frac{\epsilon}{8})\Big) , \;\;
\Big( {\gamma}_{1}(\frac{\epsilon}{8}), {\gamma}_{i2}(\frac{\epsilon}{8})\Big) , \;\;
\Big( {\gamma}_{1}(\frac{\epsilon}{8}), {\gamma}_{2}(\frac{\epsilon}{8})\Big) , \;\;
\Big( y_{01},y_{02}\Big) 
\end{displaymath}
of length
\begin{displaymath}
{\tilde{\Gamma}}^c \; < \; \frac{\epsilon}{4} \; + \; \frac{\epsilon}{8} \; + \;
\frac{\epsilon}{8} \; + \; \frac{\epsilon}{4} \; = \epsilon.
\end{displaymath}
\hfill $\Box$

\begin{corollary} \label{corollary-geodesic}
$(Y,d)$ is locally compact, complete and hence proper and geodesic.
\end{corollary}
{\bf Proof:} Since the Busemann functions $b_i$, $i=1,2$, are continuous, $Y$ is a closed subset 
of the locally compact space $(X_1\times X_2,d_m)$ and therefore locally compact itself when endowed with 
the induced metric $d_m|_Y$. Thus from Lemma \ref{lemma-same-topology} it follows that $(Y,d)$ is also locally compact. \\
Every Cauchy-sequence in $(Y,d)$ is a Cauchy-sequence in $(Y,d_m|_Y)$. But $(Y,d_m|_Y)$ is complete and therefore 
the Cauchy-sequence converges in $(Y,d_m|_Y)$. Now the proof of Lemma \ref{lemma-same-topology} yields convergence in 
$(Y,d)$. Hence $(Y,d)$ is complete. \\
Finally every locally compact, complete length space is proper and geodesic (see e.g. Proposition I.3.7. in \cite{brihae}). 
\hfill $\Box$ \\

Let $p,p'\in Y$ and let 
\begin{displaymath}
\sigma :[0,d(p,p')]\longrightarrow Y \; , \hspace{0.5cm} \sigma (t) \; = \; \Big( {\sigma}_1(t),{\sigma}_2(t)\Big)
\end{displaymath}
be a unit speed geodesic from $p$ to $p'$. We want to compare $\sigma$ with the curve ${\Gamma}(p,p')$.
To have the same domain, we modify ${\Gamma}_{pp'}$ a little: Let $a,b$ as above and 
\begin{eqnarray*}
a^* & = & a \; + \; \frac{1}{2} \Big( d(p,p') \; - \; d_m(p,p')\Big) \; , \\
b^* & = & b \; + \; \frac{1}{2} \Big( d(p,p') \; - \; d_m(p,p')\Big) .
\end{eqnarray*}
We define ${\Gamma}^*={\Gamma}^*_{pp'}
:[0,d(p,p')]\longrightarrow Y$ via
\begin{displaymath}
{\Gamma}^*_{pp'} (t) \; := \; {\gamma}|_{[0,a^*]} \; * \; {{\gamma}'|_{[0,b^*]}}^{-1}
\end{displaymath}
and prove the 
\begin{proposition} \label{prop-sigma-gamma*}
For $\sigma$ and ${\Gamma}^*$ as above it holds
\begin{displaymath}
d\Big( \sigma (t),{\Gamma}^*(t)\Big) \; \le \; 500 \delta .
\end{displaymath}
\end{proposition}

To simplify the arguments, we will assume:
\begin{description}
\item[(a)] $d(p_1,p_1')\ge d(p_2,p_2')$ which implies $a=a_1$ and $b=b_1$,
\item[(b)] $B_i(p_i)=0$, $B_i({\gamma}_i(t))=-t$, and $B_i({\gamma}_i'(t))=-t+(b-a)$.
\end{description}
We can assume this without loss of generality: $(a)$ by interchanging the factors and $(b)$ by
adding the same constant to both Busemann functions. \\

The first step in the proof of Proposition \ref{prop-sigma-gamma*} is the

\begin{lemma} \label{lemma-M}
There exists $t_0\in [0,d(p,p')]$ with $d_1({\sigma}_1(t_0),{\gamma}_1(a))\le 30\delta$.
\end{lemma}
{\bf Proof:} Consider the ideal triangle $p_1,p_1',u_1={\gamma}_1(\infty )={\gamma}_1'(\infty )$ in $X_1$ with
points $\tilde{p}_1'={\gamma}_1(a)$, $\tilde{p}_1={\gamma}'(b)$, $\tilde{u}_1={\gamma}_{p_1,p_1'}(a)\in \overline{p_1,p_1'}$
of pairwise distance $\le 8\delta$. Choose $t_0$ such that 
\begin{displaymath}
M \; := \; d\Big( {\sigma}_1(t_0),\tilde{u}_1\Big) \; = \; \min\limits_t d\Big( {\sigma}_1(t),\tilde{u}_1(t)\Big) .
\end{displaymath}
By Lemma \ref{lemma-morse-estimate} we have $L_{d_1}({\sigma}_1)\ge (a+b)+\frac{1}{20\delta}M^2$. Since
$L_{d_1}({\sigma}_1)\le L_{d_m}(\sigma )\le (a+b)+20\delta$ we obtain $M\le 20\delta$ and thus the result.
\hfill $\Box$ \\

We decompose $\sigma$ into two pieces $\sigma =\bar{\alpha} * \bar{\beta}$ where $\bar{\alpha}=\sigma |_{[0,t_0]}$ and
$\bar{\beta}=\sigma |_{[t_0,d(p,p')]}$, write $\bar{\alpha}=(\bar{\alpha}_1,\bar{\alpha}_2)$ and 
$\bar{\beta}=(\bar{\beta}_1,\bar{\beta}_2)$ and prove the 

\begin{lemma} \label{lemma-Lbaralhphai}
With the notation above it holds
\begin{displaymath}
\Big| L(\bar{\alpha}_i) \; - \; a \Big| \; \le \; 50\delta \hspace{0.5cm} \mbox{and} \hspace{0.5cm}
\Big| L(\bar{\beta}_i) \; - \; b \Big| \; \le \; 50\delta \;\;\; \mbox{for} \;\; i=1,2.
\end{displaymath}
\end{lemma}
{\bf Proof:} With $M$ as in Lemma \ref{lemma-M} we compute
\begin{displaymath}
L(\bar{\alpha}_i) \; \ge \; \Big| B_i\Big( {\sigma}_i(t_0)\Big) \; - \; B_i(p_i)\Big| \; = \;
\Big| B_1\Big( {\sigma}_1 (t_0)\Big) \Big| \; \ge \; a \; - \; M \; - 8\delta \; \ge \; a \; - \; 30\delta . 
\end{displaymath}
and 
\begin{displaymath}
L(\bar{\beta}_i) \; \ge \; \Big| B_i\Big( {\sigma}_i(t_0)\Big) \; - \; B_i(p_i')\Big| \; = \;
\Big| B_1\Big( {\sigma}_1 (t_0) \Big) \; - \; (b-a)\Big| \; \ge \; b \; - \; 30\delta . 
\end{displaymath}
Since $L(\bar{\alpha}_i) + L(\bar{\beta}_i) \le a+b+20\delta$ by Lemma \ref{lemma-estimate-Gammac} we obtain the result.
\hfill $\Box$ 

\begin{lemma} \label{lemma-M2}
With $t_0$ as in Lemma \ref{lemma-M} it holds $d_2({\sigma}_2(t_0),{\gamma}_2(a))\le 100\delta$.
\end{lemma}
{\bf Proof:} Consider the ideal triangle $p_2, q={\sigma}_2(t_0)=\bar{\alpha}_2(t_0), u_2={\gamma}_2(\infty )$ in $X_2$
with corresponding points $\tilde{q}={\gamma}_2((q\cdot B_2)_{p_2}), \tilde{u}_2, \tilde{p}_2$.
Since
\begin{displaymath}
(q\cdot B_2)_{p_2} \; = \; \frac{1}{2}
\Big(
d(p_2,\bar{\alpha}_2(t_0)) \; + \; B_2(p_2) \; - \; B_2(\bar{\alpha}_2(t_0)) 
\Big)
\end{displaymath}
and $|L(\bar{\alpha}_1-a)|\le 50\delta$, $|B_2(p_2)-B_2(\bar{\alpha}_e(t_0))|\ge a-30\delta$ by Lemma 
\ref{lemma-Lbaralhphai}, we see
\begin{displaymath}
\Big| (q\cdot B_2)_{p_2} \; - \; a\Big| \; \le \; 40\delta
\end{displaymath}
and
\begin{displaymath}
d(\tilde{u}_2,\tilde{p}_2) \; = \; \frac{1}{2}
\Big(
d(p_2,\bar{\alpha}_2(t_0)) \; + \; B_2(\bar{\alpha}_2(t_0)) \; - \; B_2(p_2)
\Big)
\; \le \; 40\delta .
\end{displaymath}
Together with Lemma \ref{lemma-xyI} we get the estimate.
\hfill $\Box$ \\

{\bf Proof of Proposition \ref{prop-sigma-gamma*}:} 
Lemmata \ref{lemma-M}, \ref{lemma-M2} and Proposition \ref{prop-d-dm} imply that $d(\sigma (t_0),\gamma (a))\le 120\delta$.
Combining some triangle inequalities we obtain $|t_0-a|\le 150\delta$ and
\begin{displaymath}
d\Big( \sigma (a^*),\gamma (a^*)\Big) \; \le 300\delta \; , \hspace{1cm}
d\Big( \sigma (a^*),\gamma (b^*)\Big) \; \le 300\delta .
\end{displaymath}
Together with Lemma $5$ we obtain Proposition \ref{prop-sigma-gamma*}.
\hfill $\Box$ 

\begin{proposition} \label{prop-yhyperbolic}
$(Y,d)$ is hyperbolic.
\end{proposition}
{\bf Proof:} By Lemma \ref{lemma-hyp-equiv2} it suffices to show that there exists a $\Delta$ such that for all 
$q,p,p'\in Y$ and all minimal geodesics $\sigma : [0,d(p,p')]\longrightarrow Y$ from $p$ to $p'$
the function $t\longmapsto d(q,\sigma (t))$ is a $\Delta$-$T$-function. By Propositions 
\ref{prop-d-dm} and \ref{prop-sigma-gamma*} it suffices to show that there exists $\Delta$ such that for all $q,p,p'$
the function
\begin{displaymath}
[a^*,b^*] \longrightarrow \mathbb{R} \hspace{1cm} t\longmapsto d_m\Big( q,{\Gamma}^*(t)\Big)
\end{displaymath}
is a $\Delta$-$T$-function, where ${\Gamma}^*={\gamma}|_{[0,a*]}*{{\gamma}'|_{[0,b^*]}}^{-1}$ as defined above. \\
Define $f,g:[a^*,b^*] \longrightarrow \mathbb{R}$ via
\begin{displaymath}
f(t) \; := \; d_1\Big( q_1,{\Gamma}^*_1(t)\Big) \hspace{0.5cm} \mbox{and} \hspace{0.5cm}
g(t) \; := \; d_2\Big( q_2,{\Gamma}^*_2(t)\Big) .
\end{displaymath}
We have to show that $\max \{ f,g\}$ is a $\Delta$-$T$-function. We use without loss of generality as above
that $d_1(p_1,p_1')\ge d_2(p_2,p_2')$, i.e. that $d_m(p,p')=d_1(p_1,p_1')$. In this case 
$d({\Gamma}_1(t),{\gamma}_{p,p'}(t))\le {\delta}'$ for all $t\in [0,a+b]$, where ${\delta}'$ only depends
on $\delta$ and not on $p,p'$. It follows that $f$ is a $\Delta$-$T$-function for some $\Delta$ only depending
on $\delta$. \\
Note further that $g|_{[0,a^*]}$ and $g|_{[a^*+\epsilon , a^*+b^*]}$ are $4\delta$-$T$-functions by Lemma 
\ref{lemma-hyp-equiv2} for every $\epsilon >0$, and hence $g|_{[a^*,a^*+b^*]}$ is a $12\delta$-$T$-function
since the jump at $a^*$ is bounded by $8\delta$. \\
Let $v=(q_2\cdot B_2)_{p_2}$, then the function 
\begin{displaymath}
[0,\infty ) \longrightarrow [0,\infty ) \; , \hspace{1cm} t\longmapsto d_2\Big( q_2,{\gamma}_2(t)\Big)
\end{displaymath}
assumes the minimum $8\delta$-close to the point $t=v$. \\
If $v>a$, hence $v\ge a^*-10\delta$, then by the properties of $\delta$-$T$-functions $g|_{[0,a^*+b^*]}$
is easily checked to be a $30\delta$-$T$-function. \\
Let us assume that $v<a$. Set $u:=d_2(q_2,p_2)-v$ and let $\tilde{\gamma}_2:[0,\infty )\longrightarrow X$ be a $B_2$-ray
starting at $q_2$. Then $d_2(\tilde{\gamma}_2(u+t),{\gamma}_2(v+t))\le 8\delta$. In particular
 $d_2(\tilde{\gamma}_2(u+(a^*-v)),{\gamma}_2(a^*))\le 8\delta$ which implies
\begin{eqnarray*}
8\delta & \ge & 
\Big|
d_2\Big( q_2,{\gamma}_2(a^*)\Big) \; - \; \Big| B_2(q_2) \; - \;  B_2\Big( {\gamma}_2(a^*)\Big) \Big| 
\Big| \\
& = & 
\Big|
d_2\Big( q_2,{\gamma}_2(a^*)\Big) \; - \; \Big| B_1(q_1) \; - \;  B_1\Big( {\gamma}_1(a^*)\Big) \Big| 
\Big| .
\end{eqnarray*}
Since $d_1(q_1,{\gamma}_1(a^*))\ge |B_1(q_1)-B_1({\gamma}_1(a))|$ we obtain
\begin{displaymath}
f(a^*) \; = \; d_1\Big( q_1, {\gamma}_1(a^*)\Big) \; \ge \;  d_2\Big( q_2, {\gamma}_2(a^*)\Big) \; - \; 8\delta \; \ge \; 
g(a^*) \; - \; 8\delta .
\end{displaymath}
Since $g|_{[0,a^*]}$ and $g|_{[a^*,a^*+b^*]}$ are $12\delta$-$T$-functions and $f$ is a ${\delta}'$-$T$-function
we see that $(g-f)_+\le 20\delta +{\delta}'$ which implies that $\max \{ f,g\}$ is a 
$20\delta +2{\delta}'$-$T$-function.
\hfill $\Box$

\begin{figure}[htbp]
\psfrag{g}{$g$}
\psfrag{f}{$f$}
\psfrag{0}{$0$}
\psfrag{v}{$v$}
\psfrag{a*}{$a^*$}
\psfrag{b*}{$b^*$}
\psfrag{le8d}{${\scriptstyle \le 8\delta}$}
\psfrag{u}{$u$}
\psfrag{v}{$v$}
\psfrag{q2}{$q_2$}
\psfrag{p2}{$p_2$}
\psfrag{g2(v)}{${\gamma}_2(v)$}
\psfrag{g2(a*)}{${\gamma}_2(a^*)$}
\psfrag{tildeg2(u)}{$\tilde{\gamma}_2(u)$}
\includegraphics[width=0.9\columnwidth]{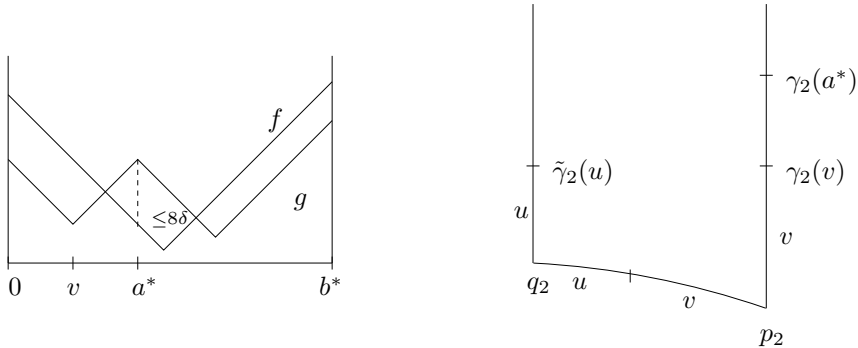}
\caption{This figure visualizes the situation in the proof of Proposition \ref{prop-yhyperbolic}.}
\end{figure}

We finally indicate how the arguments of this section have to be modified to prove Theorem \ref{theo-main}. In 
the case of Theorem \ref{theo-main} let $B_1:=d(z_1,\cdot )$, $B_2:=d(z_2,\cdot )$ and $z=(z_1,z_2)$. The $B_i$-rays correspond
to the geodesics ${\gamma}_i={\gamma}_{p_iz_i}$, ${\gamma}'_i={\gamma}_{p_i'z_i}$. Let further
\begin{eqnarray*}
\gamma : [0,d(p,z)] & \longrightarrow & Y \; , \hspace{1cm} \gamma =({\gamma}_1,{\gamma}_2) \\
{\gamma}' : [0,d(p,z)] & \longrightarrow & Y \; , \hspace{1cm} {\gamma}' =({\gamma}_1',{\gamma}_2') \; ,
\end{eqnarray*}
$a_i:=(p_i'\cdot z_i)_{p_i}$ and $b_i:=(p_i \cdot z_i)_{p_i'}$. \\
While the definition of $\Gamma$ is analog to the one in the case of Theorem \ref{theo-main2}, the definition
of ${\Gamma}^c_{pp'}$ needs to be slightly modified in the case that $a+2\delta > d(p,z)$. In that case just take
\begin{displaymath}
{\Gamma}^c_{pp'} \; := \; \gamma |_{[0,a+\tau ]} \; * \; {\gamma |_{[0,b+\tau ]}}^{-1} , 
\end{displaymath}
where $a\le \tau < 2\delta$ is chosen such that $\gamma (\tau )={\gamma}' (\tau )=z$. \\
The proof of Lemma \ref{lemma-same-topology} stays valid in the case $y_0\neq z$, $\epsilon < d(y_0,z)$. In the 
case $y_0=z$ the result is obvious.

%%%%%%%%%%%%%%%%%%%%%%%%%%%%%%%%%%%%%%%%%%%%%%%%%%%%%%%%%%%%%%%%%%%%%%%%%%%%%%%%%%%%%%%%%%%%%
%%%%%%%%%%%%%%%%%%%%%%%%%%%%%%%%%%%%%%%%%%%%%%%%%%%%%%%%%%%%%%%%%%%%%%%%%%%%%%%%%%%%%%%%%%%%%

\section{The boundary of $Y$}
\label{sec-boundary}

In the case of Theorem \ref{theo-main} it is easy to see that $\partial Y=\partial X_1\times \partial X_2$. The situation of 
Theorem \ref{theo-main2} is more interesting: \\
We study $\partial Y$ and show that it is homeomorphic to $\partial X_1 \wedge \partial X_2$. \\
Recall that the Busemann functions $B_i$ are defined as $B_i(x)=\lim_{t\rightarrow \infty}(d_i(x,{\gamma}_i(t))-t)$,
where ${\gamma}_i:[0,\infty )\longrightarrow X_i$ is a geodesic ray, $i=1,2$. \\
Let $z_i={\gamma}_i(0)$ and $u_i=[\{ {\gamma}_i(k)\}]\in \partial X_i$, $i=1,2$. Let further $z=(z_1,z_2)\in Y$ and
$\gamma (t)=({\gamma}_1(t),{\gamma}_2(t))$. Then $\gamma$ is a ray with $\gamma (0)=z$ and clearly the 
Busemann function $B:Y\longrightarrow \mathbb{R}$ of this ray is $B(y_1,y_2)=B_1(y_1)=B_2(y_2)$. \\
Let $u:=[\{ \gamma (k)\} ]\in \partial Y$. Using the results of Sections \ref{sec-preliminaries} and \ref{sec-hyp-product}
there exists a $\Delta$ such that the following holds:
\begin{description}
\item[(1)] $Y$ is $\Delta$-hyperbolic,
\item[(2)] $d({\Gamma}^*_{xy}(t),{\gamma}_{xy}(t)) \le \Delta$,
\item[(3)] $d({\Gamma}^*_{xy}(t),{\Gamma}^*_{xz}(t)) \le \Delta \;\;\; \forall \; 0\le t \le (y\cdot z)_x$,
\item[(4)] $d({\Gamma}^*_{zx}(t),{\gamma}(t)) \le \Delta \;\;\; \forall \; 0\le t\le (x\cdot B)_z$,
\item[(5)] $|d_m(x,y)-d(x,y)|\le \Delta$.
\end{description}
For any point $v\in \partial Y$ consider a geodesic ray $\sigma :[0,\infty )\longrightarrow Y$ such that $\sigma (0)=z$
and $v=[\{ \sigma (k)\} ]$. Consider the curves ${\Gamma}^*_{z\sigma (k)}$ with 
$d({\Gamma}^*_{x\sigma (k)}(t),\sigma (t))\le \Delta$ for $0\le t\le k$. By Section \ref{sec-hyp-product}
${\Gamma}^*_{z\sigma (k)}={\gamma}|_{[0,a_k^*]}*{{\gamma}^k|_{[0,b_k^*]}}^{-1}$, where 
$|a_k^*-(\sigma \cdot B)_z|\le \Delta$ and ${\gamma}_i^k:[0,\infty )\longrightarrow X_i$ is a $B_i$-ray with
${\gamma}_i^k(0)={\sigma}_i(k)$ and $d({\gamma}_i^k(b_k^*),{\gamma}_i(a_k^*))\le \Delta$. \\
We distinguish two cases: \\
{\it (Case 1) There exists a subsequence $\{ a^*_{k_j}\}$ with $\lim_{j\rightarrow \infty}$}. \\
Then $\lim_{j\rightarrow \infty}(\sigma (k_j)\cdot B)_z=\infty$ and thus
$\liminf_{j,l\rightarrow \infty}(\sigma (k_j)\cdot \gamma (l))_z=\infty$, which implies 
$\{ \sigma (k_j) \} \sim \{ \gamma (j) \}$, hence $[\{ \sigma (k_j)\} ]=u$. Since $\sim$ is an equivalence
relation and clearly $\{ \sigma (k_j) \} \sim \{ \sigma (j) \}$ we also see that $\lim_{k\rightarrow \infty}a_k=\infty$. \\
{\it (Case 2) $\{ a_k^*\}$ is bounded}. \\
Then it holds $b_k^*\longrightarrow \infty$. Reparameterize ${\gamma}_i^k:[0,\infty )\longrightarrow X_i$ as
\begin{displaymath}
\bar{\gamma}_i^k:[a_k^*-b_k^*,\infty ) \; , \hspace{1cm} \bar{\gamma}_i^k(t) \; = \; {\gamma}_i^k(t+b_k^*-a_k^*) .
\end{displaymath}
By the discussion of Sections \ref{sec-preliminaries} and \ref{sec-hyp-product}
the $\bar{\gamma}_i^k$ converge to a complete geodesic ${\gamma}_i^*$ with $B_i({\gamma}_i(t))=-t$
and $d_i({\gamma}_i(a^*_k),{\gamma}^*(a_k^*))\le 8\delta$. Clearly we have $[\{ {\gamma}_i^*(k)\} ]=u_i$. 
Let $v_i:=[\{ {\gamma}_i^*(-k)\} ]\in \partial X_i\setminus \{ u_i\}$. \\

From the discussion of Case 2 it is not difficult to show that the map
\begin{displaymath}
\begin{array}{ccc}
\partial Y \setminus \{ u\} & \longrightarrow & (\partial X_1\setminus \{ u_1\} ) \; \times \; 
(\partial X_2\setminus \{ u_2\} ) \\
& & \\
v & \longmapsto & (v_1,v_2)
\end{array}
\end{displaymath}
is a homeomorphism, which by the discussion of Case 1 extends naturally to a homeomorphism 
\begin{displaymath}
\partial Y \; \longrightarrow \; \partial X_1 \; \wedge \; \partial X_2 .
\end{displaymath}

%%%%%%%%%%%%%%%%%%%%%%%%%%%%%%%%%%%%%%%%%%%%%%%%%%%%%%%%%%%%%%%%%%%%%%%%%%%%%%%%%%%%%%%%%%%%%
%%%%%%%%%%%%%%%%%%%%%%%%%%%%%%%%%%%%%%%%%%%%%%%%%%%%%%%%%%%%%%%%%%%%%%%%%%%%%%%%%%%%%%%%%%%%%

{\footnotesize UNIVERSIT\"AT Z\"URICH, MATHEMATISCHES INSTITUT, WINTERTHURERSTRASSE 190, 
CH-8057 Z\"URICH, SWITZERLAND \\
E-mail addresses: $\;\;\;\;\;$ foertsch@math.unizh.ch $\;\;\;\;\;$ vschroed@math.unizh.ch}

%%%%%%%%%%%%%%%%%%%%%%%%%%%%%%%%%%%%%%%%%%%%%%%%%%%%%%%%%%%%%%%%%%%%%%%%%%%%%%%%
%%%%%%%%%%%%%%%%%%%%%%%%%%%%%%%%%%%%%%%%%%%%%%%%%%%%%%%%%%%%%%%%%%%%%%%%%%%%%%%%


\begin{thebibliography}{cccccccc}
%\bibitem[AlBi]{albi} S.B. Alexander \& R.L. Bishop, {\it Warped products of Hadamard spaces}, manuscripta math. 96, 1998, 487-505
\bibitem[BeKa]{beka} N.Benakli \& I. Kapovich, {\it Boundaries of hyperbolic groups}, preprint math.GR/0202286
\bibitem[BrFa]{brfa} N. Brady \& B. Farb, {\it Filling-Invariants at Infinity for Manifolds of Nonpositive Curvature},
Trans. Amer. Math. Soc., Vol. 350, Num. 8, 1998
\bibitem[BriH]{brihae}  M. Bridson \& A. Haefliger, {\it Metric spaces of non-positive curvature}, Springer Verlag Berlin 1999
\bibitem[BuBuI]{bubui} D. Burago \& Y. Burago \& S. Ivanov, {\it A course in Metric Geometry}, Graduate Studies in Mathematics, 
Vol. 33, Amer. Math. Soc., 415pp, 2001
%\bibitem[BujS]{bus} S.Bujalo \& V.Schroeder, {\it Hyperbolic rank and subexponential corank of metric spaces}, 
%Geom. funct. anal., vol. 12, 2002, 293-306
\bibitem[F1]{f1} T. Foertsch, {\it Bilipschitz Embeddings of Negative Sectional Curvature in Products of Warped Product Manifolds},
Proc. Amer. Math. Soc. 130, 2089-2096, (2002) 
\bibitem[F2]{f2} {\it The Hyperbolic Rank of Homogeneous Hadamard Manifolds}, to appear in the Manuscripta Mathematicae
\bibitem[FS]{fs} T. Foertsch \& V. Schroeder, {\it Hyperbolic Rank of Products}, preprint
%\bibitem[G1]{g1} M. Gromov, {\it Hyperbolic groups}, Essays in group theory (S.M.Gersten, ed), Springer Verlag MSRI Publ. 8 1987,
%75-263
%\bibitem[G2]{g2} {\it Asymptotic invariants of infinite groups}, in: G.A.Niblo, M.A. Roller (eds)
%Geometric group theory, vol.2.,London Math. Soc. Lect. Note series no.182, Cambridge Univ.Press 1993, 1-295
\bibitem[L1]{l1} E. Leuzinger, {\it Corank and Asymptotic Filling-Invariants for Symmetric Spaces},  Geom. Funct. Anal., 
vol.10, no4, 863-873, (2000)
\bibitem[L2]{l2} {\it Bi-Lipschitz Embeddings of Trees into Euclidean Buildings}, preprint
\bibitem[M]{m} C.R.F. Maunder {\it Algebraic Topology}, Cambridge University Press, 1980
%\bibitem[R]{r} Yu.G. Reshetnyak, {\it On the theory of spaces of curvature not greater than K}, Mat. Sbornik 52, 1960, 789-798 
\end{thebibliography}
\end{document}